
\baselineskip=14pt
\parskip=10pt

\magnification=\magstephalf

\def\1{{\overline{1}}}
\def\2{{\overline{2}}}
\parindent=0pt
\overfullrule=0in

\def\frac#1#2{{#1 \over #2}}

\centerline
{
 \bf Going Back to Neil Sloane's FIRST LOVE (OEIS Sequence A435):
}
\centerline
{
\bf On the  Total Heights in Rooted Labeled Trees
}
\rm
\bigskip
\centerline
{\it By Shalosh B. EKHAD and Doron ZEILBERGER}

\quad\quad\quad {\it Dedicated to Neil Sloane and the many contributors to the OEIS. Keep up the good work!}

{\bf Preface}

According to the ``brief history'', {\tt  http://oeis.org/wiki/Welcome\#OEIS:\_Brief\_History}, written
by Neil Sloane himself: 

{\it
``The sequence database was begun by Neil J. A. Sloane in early 1964 when he 
was a graduate student at Cornell University in Ithaca, NY. He had encountered a sequence of numbers while working on his dissertation, 
namely 1, 8, 78, 944, ... (now entry A000435 in the OEIS), and was looking for a formula for the n-th term, in order to determine 
the rate of growth of the terms.''}

That first sequence, now entry {\tt http://oeis.org/A000435} in the more than quarter-million sequences strong OEIS [Sl1], is expressible by the formula
$$
(n-1)! \sum_{k=0}^{n-2} \frac{n^k}{k!} \quad .
$$
It appears on page 119 in Sloane's Ph.D. thesis [Sl2], 
and in a joint paper with John Riordan [RS], they showed  that this is the sum of the the ``total heights'', taken over all
labeled rooted trees with $n$ vertices, divided by $n$.

{\bf Rooted Labeled Trees}

Suppose that you have a society with $n$  individuals, let's call them $1, \dots , n$, where there is a {\bf unique} ``big boss'' (the ``root'').
Every  member of the society, except the big boss, has a unique {\it immediate supervisor}. Some people (``leaves'')
have no one reporting to them, but the set of immediate subordinates of each supervisor is {\it unordered}, i.e.
they are considered of equal status in the ``pecking order''. Of course, no one can be their own (immediate or indirect)
supervisor. How many such hierarchies are possible?

If you draw the hierarchy with a directed edge between any member and his immediate supervisor, you would get
a {\bf labeled rooted tree}. Arthur Cayley[C] famously proved that the number of labeled trees on $n$ vertices
is $n^{n-2}$, hence the number of rooted labeled trees is $n \cdot n^{n-2}=n^{n-1}$.

There are many proofs of this result, the nicest one is due to Andr\'e Joyal[J] (see also [LZ]). Another one is using
{\bf Lagrange Inversion} (see [Z1] for a nice exposition), and that's the one needed for the
present article. Let's review it.

Let  $r(n)$ be the number of labeled rooted trees with $n$ vertices, and consider the {\bf exponential generating function}
$$
R(x) \, := \, \sum_{n=0}^{\infty} \frac{r(n)}{n!} x^n  \quad .
$$

If the degree of the root is $k$, then deleting it gives us a {\bf set} (i.e. unordered collection) of
smaller rooted labeled trees (with disjoint labels), 
that by general generatingfunctionology has exponential generating function
$x \frac{R(x)^k}{k!}$ (we divide by $k!$ since the $k$ subtrees are unordered). Summing over all
possible $k \geq 0$ , we get
$$
R(x) = x \sum_{k=0}^{\infty}   \frac{R(x)^k}{k!} \, = \,    x e^{R(x)} \quad .
$$
We have just established a {\bf functional equation} for the formal power series, $R(x)$:
$$
R(x)= x e^{R(x)} \quad .
$$

For any formal power series $f(t)$, let $[t^n] f(t)$  denote the coefficient of $t^n$ in $f(t)$.

Recall the versatile

{\bf Lagrange Inversion Theorem}:
If $R(x)$ and $\Phi(z)$ are formal power series, starting at $x$ and
$z^0$ respectively, then $R(x)=x \Phi(R(x))$ implies
$[x^n] R(x) = \frac{1}{n} [z^{n-1}] \Phi(z)^n$ \quad .

In our case $\Phi(z)=e^z$, hence $\Phi(z)^n=e^{nz}$, whose coefficient of $z^{n-1}$ is $\frac{n^{n-1}}{(n-1)!}$.
Hence the coefficient of $x^n$ in $R(x)$, $r(n)/n!$, equals
$$
\frac{r(n)}{n!} = \frac{1}{n} \cdot \frac{n^{n-1}}{(n-1)!} \quad,
$$
entailing that, indeed,
$$
r(n)=n^{n-1} \quad .
$$

We have just proved that there are exactly $n^{n-1}$ labeled rooted trees with $n$ vertices.

So much for {\it naive} counting, but there is a lot of diversity among these hierarchies.
One extreme is that everyone, except the big boss, reports directly to the big boss, so the
``distance'' to the root is always $1$ and the sum of the distances is $n-1$.
There are only $n$ such trees, since once you have chosen the root (the ``big boss'')
there is nothing to do. This is the most democratic rooted tree.

The other extreme is that the  hierarchy  is {\bf totally} ordered. Every vertex has only one
subordinate, except the one at the very bottom, that has none. Now the sum of the distances to the root is
$0+1+2+ \dots +(n-1)= n(n-1)/2$ (and hence  the average distance is $(n-1)/2$), and there are
$n!$ such trees. Such trees are the most authoritarian, there is a clear ranking, and no one
is of equal status.

Hence a natural measure of how ``authoritarian'' a rooted tree is, is the sum of the heights (distances to the root)
taken over all vertices. Let's  define the {\bf weight-enumerator} of the set of labeled rooted trees on
$n$ vertices by
$$
J_n(y) :=\sum_{T} y^{TotalHeight(T)} \quad ,
$$
where the sum is taken over the set of rooted labeled trees on $n$ vertices.
Of course $J_n(1)=n^{n-1}$, but can we we find an explicit expression for $J_n(y)$ in terms of $n$ and $y$?
Probably not! Still it would be nice to have an efficient algorithm to generated as many terms of the
polynomial sequence $\{J_n(y) \}$ as possible, and also be able to find an explicit expression for
$J_n'(1)$, since the important quantity ``average total height'' is given by $J_n'(1)/n^{n-1}$.
In fact that was Neil Sloane's original  motivation, that lead to sequence A435.

{\bf Weighted Counting According to Total Height}

Riordan and Sloane[RS] define the formal power series of the {\bf two} variables $x$ and $y$
$$
J(x,y)= \sum_{n=1}^{\infty} J_n(y) \frac{x^n}{n!} \quad ,
$$
that is the exponential generating function of the sequence of polynomials $\{ J_n(y) \}$. Of course
$J(x,1)=R(x)$.

Using the same generatingfunctionology argument, it is not hard to show (as done in [RS]) that
$J(x,y)$ satisfies the {\bf functional equation}
$$
J(x,y)=x e^{J(xy,y)} \quad .
\eqno(FE)
$$
Alas, now Lagrange Inversion is no longer applicable, and there is no way to recover $J_n(y)$ explicitly.

But what about  $J_y(x,1)$?  (i.e. $\frac{\partial}{\partial y} J(x,y)$ evaluated at $y=1$).

Let's  differentiate Eq. $(FE)$ with respect to $y$, 
recalling the chain rule from multivariable calculus. We get
$$
J_y(x,y) =x e^{J(xy,y)} \cdot \frac{\partial}{\partial y} J(xy,y)=
J(x,y) \cdot \frac{\partial}{\partial y} J(xy,y)=
J(x,y) \cdot \left ( \frac{\partial(xy) }{\partial y} \cdot \frac{\partial}{\partial (xy)} J(xy,y) + \frac{\partial}{\partial y} J(xy,y) \right )
$$
$$
=J(x,y) \cdot \left ( x J_x(xy,y) + J_y(xy,y) \right ) \quad .
$$
Now plug-in $y=1$ to get
$$
J_y(x,1)=J(x,1) (x J_x(x,1)+ J_y(x,1)) \quad .
$$

But $J(x,1)$ is what we called above $R(x)$, and $J_x(x,1)$ is $R'(x)$, hence
$$
J_y(x,1)= x R(x) R'(x) + R(x) J_y(x,1) \quad .
$$
Solving for $J_y(x,1)$ we get
$$
J_y(x,1)= \frac{xR(x)R'(x)}{1-R(x)} \quad .
$$

It would be nice to express $J_y(x,1)$ in terms of $R(x)$ only, but this is easy.

Differentiating the functional equation $R(x)=x \, e^{R(x)}$ with respect to $x$,
we get, by the product rule and chain rule (this time Calculus I suffices)
$$
R'(x)= e^{R(x)} + x e^{R(x)} R'(x)=\frac{R(x)}{x} + R(x) R'(x) \quad .
$$
Solving for  $R'(x)$ we get
$$
R'(x)= \frac{R(x)}{x(1-R(x))} \quad .
$$

We note, for the future, that by repeated differentiation (using the quotient rule and the chain rule and repeatedly using
that very same equation $R'(x)= \frac{R(x)}{x(1-R(x))}$)
enables us to express any derivative of $R(x)$, $R^{(j)}(x)$, as rational function of $R(x)$ and $x$
with denominator that has the form $(1-R(x))^{2j-1}$.

Substituting $R'(x)= \frac{R(x)}{x(1-R(x))}$ into $J_y(x,1)= \frac{xR(x)R'(x)}{1-R(x)}$ gives
$$
J_y(x,1)= \frac{R(x)^2}{(1-R(x))^2} \quad .
$$

Now it is time to invoke (see, e.g., [Z1])

{\bf The Generalized Lagrange Inversion Theorem:}
If $u(t)$ and $\Phi(z)$ are formal power series starting at $t$ and
$z^0$ respectively, and $G(z)$ is yet another formal power series,
then $u(t)=t \Phi(u(t))$ implies
$[t^n] G(u(t)) = (1/n) [z^{n-1}] \, G'(z) \Phi(z)^n$ \quad .

Here $G(z)=\frac{z^2}{(1-z)^2}$ and
hence $G'(z)=2\,{\frac {z}{ \left( 1-z \right) ^{3}}}\,\,$.

Hence $J_n'(1)/n!$, the coefficient of $x^n$ in  $\frac{R(x)^2}{(1-R(x))^2}$, is $1/n$ times the coefficient of
$z^{n-1}$ in
$$
{\frac {2z}{ \left( 1-z \right) ^{3}}} \cdot e^{nz}  \quad,
$$

which is the coefficient of $z^n$ in
$$
{\frac {2z^2}{ \left( 1-z \right) ^{3}}} \cdot e^{nz}  \quad .
$$

But
$$
{\frac {2z^2}{ \left( 1-z \right) ^{3}}} = \sum_{k=0}^{\infty} (k-1)k \, z^k \quad ,
$$
hence $\frac{J_n'(1)}{n!}$ is $\frac{1}{n}$ times the coefficient of $z^n$ in the formal power series
$$
\left ( \sum_{k=0}^{\infty} (k-1)k \, z^k \right ) \cdot \left ( \sum_{s=0}^{\infty} \frac{n^s}{s!} \right ) \quad .
$$
Hence
$$
\frac{J_n'(1)}{n!} \, = \, \frac{1}{n} \sum_{k=0}^{n-2} \frac{(n-k-1)(n-k) \, n^{k}}{k!} 
\quad .
$$
Noting that
$$
(n-k)(n-k-1) \, = \,n \,( n-1) \,- \, 2( n-1) \, k + k\,(k-1) 
\quad ,
$$
simple routine algebra leads to
$$
J_n'(1)= n! \sum_{k=0}^{n-2} \frac{n^k}{k!} \quad .
$$
Hence the average total height among all labeled rooted trees with $n$ vertices is
$\frac{n!}{n^{n-1}} \sum_{k=0}^{n-2} \frac{n^k}{k!}$, that we will call $W_n$ (please be warned
that our notation differs from that of [RS], their $W_n$ is $n^{n-1}$ times our $W_n$).

As noted in [RS], $W_n$ is asymptotic (thanks to Ramanujan and Watson, see [W]) to $n^{3/2}\sqrt{\pi/2}$.
We have just reproved, in much more detail than in [RS] (and a somewhat different proof):

{\bf Theorem 1} (Riordan-Sloane [RS])
The average total height among all rooted labeled trees on $n$ vertices equals
$\frac{n!}{n^{n-1}} \sum_{k=0}^{n-2} \frac{n!}{k!}$ and is asymptotically  $n^{\frac{3}{2}}\sqrt{\pi/2}$.

{\bf Enter Computers}

So much can be done by mere humans, but the average is only the {\bf most basic} statistical information about a random variable.
What about the variance? (and hence ``coefficient of variation'') skewness? kurtosis? and  higher moments?
Is there a limiting scaled distribution?

In order to find explicit expressions for higher moments, we need to first find higher factorial moments. The $r$-th factorial moment is
$J_n^{(r)}(1)$, and once we know the first $r$ factorial moments we can, by standard theory (see [Z2]), get the
{\it moments}, and from them, easily, the {\it moments-about-the-mean}.

But how can we do that? It turns out that the same method that we described above still works, but very soon gets
very tedious for humans. To get the second factorial moment, we have to differentiate $(FE)$ twice, plug-in $y=1$ and
get an expression for $J_{yy}(x,1)$
(that is the exponential generating function of $J_n''(1)$,) in terms of $R(x),R'(x)$ and $R''(x)$. We already noted that
each derivative of $R(x)$ can be expressed as rational function of $R(x)$,
so at the end it can be expressed in terms of $R(x)$ alone, and we can use the Generalized Lagrange Inversion
Formula, as we did above.

Maple knows the chain rule for multi-variable functions, so all this can be done automatically and seamlessly.
Also one can teach Maple how to use  generalized Lagrange Inversion, and perform all the steps.

All this is implemented in the Maple package {\tt A435.txt} available from the front of this article
{\tt http://www.math.rutgers.edu/\~{}zeilberg/mamarim/mamarimhtml/a435.html}\quad .

{\bf Some Computer-Generated Theorems}

In the theorems below
$$
W_n \, := \, \frac{n!}{n^{n-1}} \sum_{k=0}^{n-2} \frac{n^k}{k!} \quad .
$$
Recall that Riordan and Sloane showed (and we reproved above) that the average total height among labeled rooted trees with $n$
vertices is $W_n$.

It follows from our algorithm that {\it every} moment can always be expressed as {\it some} polynomial in $n$ and $W_n$,
but they get more and more complicated for higher moments. Below we stae rigorously-proved explicit expressions for
the first four moments, as well as the implied asymptotics and the limits of the $\alpha$-coefficients, i.e.
the limits of the standardized moments. More moments can be found in the output files mentioned later.

{\bf Theorem 2.}
 The variance of the random variable ``total height'' on the set of rooted labeled trees on $n$ vertices is given explicitly by
$$
 -{W_{{n}}}^{2} -{\frac {17}{6}}\,nW_{{n}} +  \frac{5}{3} \, {n}^{2}(n-1) \quad ,
$$
and its asymptotics is $(\frac{5}{3} -\frac{\pi}{2})n^2$. Hence the limit of the coefficient of variation
(the mean over the standard-deviation), as $n$ goes to infinity, is
$$
\frac{\sqrt{2}}{6}\,{\frac {\sqrt {-18\,\pi +60}}{\sqrt {\pi }}} \, = \, 0.2470484847\dots \quad .
$$
Note in particular that there is no ``concentration about the mean''.

{\bf Theorem 3.}
 The third moment about the mean of the random variable ``total height'' on the set of rooted labeled trees on $n$ vertices is given explicitly by
$$
2\,{W_{{n}}}^{3}+\frac{17}{2} \,n{W_{{n}}}^{2}+ 
\left( -{\frac {25}{8}}\,{n}^{3} + {\frac {277}{24}}\,{n}^{2} - {\frac {1}{60}}\,n \right) W_{{n}}
-{\frac {151}{30}}\,{n}^{4}
+ {\frac {76}{15}}\,{n}^{3}
- {\frac{1}{30}} \,{n}^{2} \quad ,
$$
and its asymptotic expression is
$$
\left( \frac{1}{2}\,\sqrt {2}{\pi }^{3/2}-{\frac {25}{16}}\,\sqrt {2}\sqrt {\pi } \right) {n}^{9/2} \quad,
$$
that is approximately $0.020795808\,{n}^{9/2}$. It follows that the limit of the {\it skewness}, as $n$ goes to infinity, is
$$
{\frac { \left( 6\,\pi -{\frac {75}{4}} \right) \sqrt {3}\sqrt {{\frac {\pi }{10-3\,\pi }}}}{10-3\,
\pi }} \, = \,  .7005665208\dots \quad .
$$
In particular we know that the limiting distribution, {\it whatever it is},  is {\bf not} normal. So ``total height'' defined on
rooted labeled trees is {\bf not} asymptotically normal.

{\bf Theorem 4.}
The fourth moment about the mean of the random variable ``total height'' on the set of rooted labeled trees on $n$ vertices is given explicitly by
$$
-3\,{W_{{n}}}^{4}-17\,n\,{W_{{n}}}^{3}+ \left( 
 \frac{5}{2} \,{n}^{3}
-{\frac {217}{6}}\,{n}^{2}+ \frac{1}{15} \,n \right) {W_{{n}}}^{2}+ \left(  {\frac {649}{80}}\,{n}^{4} -{\frac {74381}{2160}}\,{n}^{3}
+{\frac {433}{2520}}\,{n}^{2}+{\frac {1}{105}}\,n \right) W_{{n}}
$$
$$
+{\frac {221}{63}}\,{n}^{6}+ {\frac {4693}{540}}\,{n}^{5} -{\frac {4651}{378}}\,{n}^{4}+{\frac {109}{1260}}\,{n}^{3}
+{\frac {2}{105}}\,{n}^{2} \quad ,
$$

and its asymptotic expression is
$$
\left( - \frac{3}{4} \,{\pi }^{2 \,   } \,       + \, \frac{5}{4} \,\pi + {\frac {221}{63}} \right) {n}^{6}  \quad ,
$$ 
that is approximately $0.032724023 n^6$. It follows that the limit of the {\it kurtosis} as $n$ goes to infinity is
$$
\frac{1}{7} \,{\frac {-189\,{\pi }^{2}+315\,\pi +884}{ \left( 10-3\,\pi  \right) ^{2}}} \,=\,  3.560394751\dots \quad ,
$$
hence the limiting distribution is  {\it leptokurtic}.

For theorems about the $5^{th}$ through the $12^{th}$ moments we refer the reader to the computer-generated article

{\tt http://www.math.rutgers.edu/\~{}zeilberg/tokhniot/oA435a12.txt} \quad .

Let us conclude by stating the limits of the scaled moments, $\alpha_k$,  for $3 \leq k \leq 9$.
$$
\alpha_3 \, = \,
{\frac{ \left( 6\,\pi -{\frac {75}{4}} \right) \sqrt {3}\sqrt {{\frac {\pi }{10-3\,\pi }}}}{10-3\,\pi }}
\, =\,
.7005665208 \dots \quad , 
$$
$$
\alpha_4 \, = \, 
{\frac{-189\,{\pi }^{2}+315\,\pi +884}{ 7 \, \left( 10-3\,\pi  \right) ^{2}}}
\, = \,
3.560394751\dots \quad , 
$$
$$
\alpha_5 \, \, = \, 
{\frac { \left( 36\,{\pi }^{2}+{\frac {75}{2}}\,\pi -{\frac {105845}{224}} \right) \sqrt {3}\sqrt {{\frac {\pi }{10-3\,\pi }}}}{ \left( 10-3\,\pi 
 \right) ^{2}}}
\, =\, 7.256376376\dots \quad , 
$$
$$
\alpha_6 \, \, = \, 
{\frac{15}{16016}}\,{\frac {-144144\,{\pi }^{3}-720720\,{\pi }^{2}+3013725\,\pi +2120320}{ \left( 10-3\,\pi  \right) ^{3}}}
\, = \, 27.68549546\dots \quad , 
$$
$$
\alpha_7 \, \, = \, 
{\frac{ \left( 162\,{\pi }^{3}+{\frac {6615}{4}}\,{\pi }^{2}-{\frac {103965}{32}}\,\pi -{\frac {101897475}{9152}} \right) \sqrt {3}\sqrt {{\frac {\pi }{
10-3\,\pi }}}}{ \left( 10-3\,\pi  \right) ^{3}}}
\, = \, 90.01702180\dots \quad, 
$$
$$
\alpha_8 \, \, = \, 
{\frac{3}{2586584}}\,{\frac {-488864376\,{\pi }^{4}-8147739600\,{\pi }^{3}-455885430\,{\pi }^{2}+86568885375\,\pi +32820007040}{ \left( 10-3\,\pi 
 \right) ^{4}}} 
$$
$$
\, = \, 358.8086679\dots \quad ,
$$
$$
\alpha_9 \, \, = \, 
{\frac{ \left( 648\,{\pi }^{4}+15795\,{\pi }^{3}+{\frac {591867}{16}}\,{\pi }^{2}-{\frac {461286225}{2288}}\,\pi -{\frac {188411947088175}{662165504}}
 \right) \sqrt {3}\sqrt {{\frac {\pi }{10-3\,\pi }}}}{ \left( 10-3\,\pi  \right) ^{4}}}
\, = \, 1460.710269\dots \quad .
$$
For the exact expressions for $\alpha_{10},\alpha_{11},\alpha_{12}$ see the output file 

{\tt http://www.math.rutgers.edu/\~{}zeilberg/tokhniot/oA435a12.txt} \quad .

Here are their floating-point approximations:
$$
\alpha_{10}  \, = \, 6498.233818 \dots \quad , 
$$
$$
\alpha_{11} \, = \, 30389.98955\dots \quad , 
$$
$$
\alpha_{12}  \, = \, 150516.4157\dots  \quad .
$$

One of us (DZ) is pledging a donation of one hundred US dollars to the OEIS Foundation, in honor of the first solver(s), for a solution to the following
challenge.

{\bf Challenge}: What is the probability density function of the  limiting scaled distribution,
as $n \rightarrow \infty$, of the random variable ``total height'' defined on the set of labeled rooted trees on $n$ vertices?

To get a glimpse of how it is supposed to look like, see the plots here:

{\tt http://www.math.rutgers.edu/\~{}zeilberg/tokhniot/oA435c.html} \quad .

\bigskip
\hrule

{\bf References}

[C] Arthur Cayley, {\it A theorem on trees}, Quart. J. Math {\bf 23}(1889), 376-378.

[J] Andr\'e Joyal, {\it  Une theorie combinatoire des s\'eries formelles}, Advances in Mathematics  {\bf 42}(1981), 1-81.

[LZ] Gyu Eun Lee and Doron Zeilberger, {\it Joyal's Proof of Cayley's Formula}, 
The Personal Journal of Shalosh B. Ekhad and Doron Zeilberger, July 16, 2012. Available on-line from  \hfill\break
{\tt http://www.math.rutgers.edu/\~{}zeilberg/mamarim/mamarimhtml/JoyalCayley.html} \quad .

[RS] John Riordan and Neil J. A. Sloane, {\it The enumeration of rooted trees by total height},
J. Australian Math. Soc. {\bf 10}(1969),  278-282. Available on line from: \hfill\break 
{\tt http://neilsloane.com/doc/riordan-enum-trees-by-height.pdf} \quad .

[Sl1] Neil J.A. Sloane, ``The On-Line Encyclopedia of Integer Sequences'', {\tt http://www.oeis.org}.

[Sl2] Neil J. A. Sloane, 
{\it ``Lengths of Cycle Times in Random Neural Networks''}, Ph.D. Dissertation, 
Cornell University, February 1967; 
also Report No. 10, Cognitive Systems Research Program, Cornell University, 1967. Available on-line from \hfill\break
{\tt http://neilsloane.com/doc/a000435\_OCR.pdf \quad } .

[W] George N. Watson, {\it Theorems stated by Ramanujan(V): Approximations connected with $e^x$}, Proc. London Math. Soc. (2), {\bf 29} (1929), 293-308.

[Z1] Doron Zeilberger, {\it Lagrange Inversion Without Tears (Analysis) (based on Henrici)},
The Personal Journal of Shalosh B. Ekhad and Doron Zeilberger, 2002. Available on-line from \hfill\break
{\tt http://www.math.rutgers.edu/\~{}zeilberg/mamarim/mamarimhtml/lag.html } \quad .

[Z2] Doron Zeilberger, {\it The Automatic Central Limit Theorems Generator (and Much More!) },
in: ``Advances in Combinatorial Mathematics: Proceedings of the Waterloo Workshop in Computer Algebra 2008 in honor of Georgy P. Egorychev'', 
chapter 8, pp. 165-174, (I.Kotsireas, E.Zima, eds., Springer Verlag, 2009.) Available on-line from \hfill\break
{\tt http://www.math.rutgers.edu/\~{}zeilberg/mamarim/mamarimhtml/georgy.html} \quad .

\bigskip
\bigskip
\hrule
\bigskip
Doron Zeilberger, Department of Mathematics, Rutgers University (New Brunswick), Hill Center-Busch Campus, 110 Frelinghuysen
Rd., Piscataway, NJ 08854-8019, USA. \hfill \break
DoronZeil at gmail dot com  \quad ;  \quad {\tt http://www.math.rutgers.edu/\~{}zeilberg/} \quad .
\bigskip
\hrule
\bigskip
Shalosh B. Ekhad, c/o D. Zeilberger, Department of Mathematics, Rutgers University (New Brunswick), Hill Center-Busch Campus, 110 Frelinghuysen
Rd., Piscataway, NJ 08854-8019, USA.
\bigskip
\hrule

\bigskip
Exclusively published in The Personal Journal of Shalosh B. Ekhad and Doron Zeilberger  \hfill \break
({ \tt http://www.math.rutgers.edu/\~{}zeilberg/pj.html})
and {\tt arxiv.org} \quad . 
\bigskip
\hrule
\bigskip
{\bf  Written:  July 19, 2016} ; 

\end